
\input amstex

\input amsppt.sty

\mag=\magstep1

\topmatter
\title General hyperplane sections of nonsingular flops in dimension 3
\endtitle
\rightheadtext {nonsingular flops}
\author Yujiro Kawamata \endauthor
\address Department of Mathematical Sciences, University of Tokyo, Hongo,
Bunkyo, Tokyo, 113, Japan \endaddress
\email kawamata\@math.s.u-tokyo.ac.jp\endemail
\endtopmatter

\document

Let $X$ be a 3-dimensional complex manifold,
and $f: X \to Y$ a proper bimeromorphic morphism to a normal complex space
which contracts an irreducible curve $C \subset X$ to a singular
point $Q \in Y$ while inducing an isomorphism $X \setminus C \simeq
Y \setminus \{Q\}$.
We assume that the intersection number with the canonical divisor
$(K_X \cdot C)$ is zero.
In this case, it is known that the singularity of $Y$ is Gorenstein terminal,
and there exists a {\it flop} $f^{\#}: X^{\#} \to Y$ ([R]), which we call
a {\it nonsingular flop} because $X$ is nonsingular.

In order to investigate $f$ analytically, we replace $Y$ by its germ at $Q$,
and consider a general hyperplane section $H$ of $Y$ through $Q$.
Then $H$ has only a rational double point,
its pull-back $L \subset X$ by $f$ is normal,
and the induced morphism $f_H: L \to H$ factors the minimal resolution
$g: M \to H$ ([R]).
The dual graph $\Gamma$ of the exceptional curves of $g$ is a
Dynkin diagram of type $A_n$, $D_n$ or $E_n$.
Let $F = \sum _{k=1}^n m_kC_k$ be the {\it fundamental cycle} for $g$ on $M$.
The natural morphism $h: M \to L$ is obtained by contracting all
the exceptional curves of $g$ except the strict transform $C_{k_0}$ of $C$.

Koll\'ar defined an invariant of $f$ called the {\it length} as the length
of the scheme theoretic fiber $f^{-1}(Q)$ at the generic point of $C$.
It coincides with the multiplicity $m_{k_0}$ of the fundamental cycle
at $C_{k_0}$.

Katz and Morrison proved the following theorem ([KM, Main Theorem]).
The purpose of this paper is to give its simple geometric proof.

\proclaim{Theorem}  Let $f: X \to Y$ be as above.
Then the singularity of the general hyperplane section $H$
and the partial resolution $f_H: L \to H$ are determined
by the length $\ell$ of $f$.
More precisely,
$H$ has a rational double point of type $A_1, D_4, E_6, E_7, E_8$, or $E_8$,
if $\ell = 1, 2, 3, 4, 5$ or $6$, respectively.
\endproclaim

We note that there is only one irreducible component of $g^{-1}(Q)$
whose multiplicity in $F$ coincides with $\ell$ in the above cases.

\demo{Proof}
Let $H'$ be another general hyperplane section of $Y$ through $Q$,
and $f_{H'}: L' \to H'$
the induced morphism.  $H$ and $H'$ have the same type of singularities, and
so do $L$ and $L'$.
Let $P_i$ and $P'_i$  be the singular points of $L$ of $L'$, respectively.

Let $D$ be the effective Cartier divisor on $L$ given by $L' \cap L$.
Then $D$ is a general member of the linear system of
effective Cartier divisors on $L$ which contain $C$ and
such that $(D \cdot C) = 0$.
In fact, if $s_0$ is the global section of
$\Cal O_L(- C) \subset \Cal O_L$ corresponding to $D$,
then from an exact sequence
$$
0 \to \Cal O_X(- L) \to \Cal O_X \to \Cal O_L \to 0
$$
there exists a lifting $s \in H^0(X, \Cal O_X)$ of $s_0$ which defines $L'$,
because $H^1(X, \Cal O_X(- L)) \mathbreak = 0$.

Let $\tilde D$ be the total transform of $D$ on $M$.
Then we can write $\tilde D = F + D'$ for some $D'$
which is reduced, nonsingular and has no common irreducible
components with $F$.
If $\Gamma$ is of type $A_n$, then $D'$ has 2 irreducible components
each of which intersects transversally one of the end components of $F$.
Otherwise, $D'$ is irreducible, and intersects transversally
a component $C_{k_1}$ of $F$ such that $m_{k_1} = 2$ ($k_1$
may be equal to $k_0$).

Let $t$ be the global section of $\Cal O_X$ corresponding to $L$.
Then $s + ct$ is also a lifting of $s_0$ for any $c \in \Bbb C$.
Let $L'(c)$ be the corresponding divisor on $X$.

Let $P$ be a point on $C$ which is different from the $P_i$.
For local analytic coordinates $\{x, y, t\}$, we can write
$s + ct = F(x, y) + t(G(x, y, t) + c)$.
For a general choice of $c$,
$(G(x, y, t) + c)\vert _C$ does not vanish at the singular points $P'_i$
of $L' = L'(0)$ other than the $P_i$,
and has only simple zeroes at some points $P''_j$.
Then $L'' = L'(c)$ has singularities only at the $P''_j$,
besides possibly at the $P_i$,
with equations of the type $x^{\ell} + ty = 0$.

If we replace $L'$ by $L''$, we conclude that $L'$ has only singularities
of type $A_{\ell - 1}$ outside the $P_i$.  We shall investigate
the singularities of $L'$ at the $P_i$ case by case.
\vskip 1pc

Let $\Gamma _i$ be the dual graph of the exceptional curves of $h$ over
$P_i$, and $F_i$ the corresponding fundamental cycle.
{}From the description of $\tilde D$ above,
we can calculate the multiplicity $d_i$
of $D$ at the point $P_i$ by $d_i = ((m_{k_0}C_{k_0} + D') \cdot F_i)$.

If $\ell = 1$ or $2$, then we can check that $d_i \le 3$,
hence $L'$ is nonsingular at the $P_i$.
Then it follows that $\Gamma = A_1$ or $D_4$, respectively.

But if $\ell \ge 3$, then $d_i$ can be bigger,
and we should look at the singularity
of $L'$ more closely.

We assume first that $\ell = 3$.
If $\Gamma = E_6$, then there is nothing to prove.
We have to prove that $\Gamma \ne E_7, E_8$.
If $\Gamma = E_7$, then $L$ has two singular points
$P_1$ and $P_2$, where $F_1$ meets $D'$.
We have 2 cases; $\Gamma _1 = A_1$ and $\Gamma _2 = A_5$,
or $\Gamma _1 = A_4$ and $\Gamma _2 = A_2$.
In the former case, $L'$ has at most $A_1$ singularity at $P_1$
because of the symmetry of $L$ and $L'$,
while being nonsingular at $P_2$, since $d_2 = 3$.
Therefore, $L'$ has simpler singularities than $L$, a contradiction.
In the latter case, it has at most $A_2$ at $P_2$.  We shall prove that
$L'$ has $A_1$ at $P_1$.

Let $\mu: X^{(1)} \to X$ be the blowing-up at $P_1$,
$E \simeq \Bbb P^2$ the exceptional divisor,
and $L^{(1)}$ (resp. $L^{(1)\prime}$) the strict transform of $L$ (resp. $L'$).
$B = L^{(1)} \cap E$ consists of 2 lines $B_1$ and $B_2$, which correspond
to the 2 end components of $F_1$.  Since their multiplicities in $F$ are 2
and $L^{(1)} \cdot E = B$,
we deduce that $\mu^*L' = L^{(1)\prime} + 2E$, and
neither of the $B_i$ are contained in $B' = L^{(1)\prime} \cap E$.
Thus the intersection of 2 conics $B$ and $B'$ is equal to
$(L^{(1)} \cap L^{(1)\prime}) \cap E$.
We see from the description of $\tilde D$
that it consists of 2 points, one at
$B_1 \cap B_2$ and the other on one component $B_1$.
Then $B'$ must be a nonsingular conic,
and $L'$ has $A_1$ singularity at $P_1$.

If $\Gamma = E_8$, then we have again 2 cases;
$\Gamma _1 = A_1$ and $\Gamma _2 = E_6$,
or $\Gamma _1 = A_7$.
In the former case, $L'$ has at most $A_1$ singularity at $P_1$, while
being nonsingular at $P_2$, since $d_2 = 3$.
In the latter case, it has $A_1$ at $P_1$ as in the case of $E_7$.
\vskip 1pc

Next we assume that $\ell = 4$.
If $\Gamma = E_7$, then there is nothing to prove.
If $\Gamma = E_8$, then we have 2 cases;
$\Gamma _1 = D_5$ and $\Gamma _2 = A_2$,
or $\Gamma _1 = A_6$ and $\Gamma _2 = A_1$.

In the former case, $L'$ has at most $A_2$ singularity at $P_2$.
By the symmetry of $L$ and $L'$, $L'$ has $D_5$ at $P_1$.
Let $\mu: X^{(1)} \to X$, $E$, $L^{(1)}$ and $L^{(1)\prime}$ as before.
$B = L^{(1)} \cap E$ is a line, and $L^{(1)} \cdot E = 2B$.
Since the corresponding curve has
multiplicity 4 in $F$, we have $\mu^*L' = L^{(1)\prime} + 2E$,
and $B$ is not contained in $B' = L^{(1)\prime} \cap E$.
$L^{(1)}$ has 2 singular points $P_1^{(1)}$ and $P_2^{(1)}$ on $B$
which are of types $A_3$ and $A_1$, respectively.
We have $B \cap B' = P_1^{(1)}$ by the description of $\tilde D$.

Let $\nu: X^{(2)} \to X^{(1)}$ be the blowing-up at $P_1^{(1)}$,
$E^{(1)} \simeq \Bbb P^2$ the exceptional divisor,
and $L^{(2)}$ (resp. $L^{(2)\prime}$) the strict transform of $L^{(1)}$
(resp. $L^{(1)\prime}$).
$B^{(1)} = L^{(2)} \cap E^{(1)}$ consists of 2 lines, and
one of the corresponding curves on $M$ has
multiplicity 3 in $F$, hence $\nu^*L^{(1)\prime} = L^{(2)\prime} + E^{(1)}$,
and $L^{(1)\prime}$ is nonsingular at $P_1^{(1)}$.
But this contradicts the symmetry of $L$ and $L'$.

In the latter case, $L'$ has at most $A_1$ singularity at $P_2$.
Let $\mu: X^{(1)} \to X$, etc., as before.
$B = L^{(1)} \cap E$ consists of 2 lines $B_1$ and $B_2$, which correspond
to the 2 end components of $F_1$.  Since their multiplicities in $F$ are
3 and 2,
$B_1$ is contained in $B' = L^{(1)\prime} \cap E$, while $B_2$ is not.
Thus we have $B' = B_1 + B'_2$ with $B_2 \ne B'_2$.
The strict transform of $C_{k_0}$ passes through the point $B_1 \cap B_2$,
a contradiction to the symmetry.

Finally, if $\ell \ge 5$, the assertion of the theorem is clear.
Q.E.D.
\enddemo

\enddocument